\documentclass{article}
\usepackage{authblk}
\usepackage{amsmath,amssymb,amsfonts}
\usepackage{graphicx}
\usepackage{textcomp}
\usepackage{float}
\usepackage[margin=1in]{geometry}
\usepackage{commath,mathtools}
\usepackage{algorithm}
\usepackage{algorithmic}
\usepackage{subfigure}
\usepackage{bigints}
\usepackage[english]{babel}
\usepackage[autostyle]{csquotes}
\usepackage{hyperref}
\usepackage{color}
\usepackage{rotating,booktabs}
\usepackage{lscape}
\usepackage[sort,nocompress]{cite}
\usepackage{enumitem}
\usepackage{multirow}
\usepackage{ulem} 
\usepackage{tikz}
\usepackage{comment}
\usetikzlibrary{shapes,arrows,chains}
\usepackage{amsthm}
\usepackage{amssymb}
\newtheorem{theorem}{Theorem}[section]

\usepackage{pifont}
\usepackage[pagewise,mathlines]{lineno}
\DeclareFontFamily{U}{matha}{\hyphenchar\font45}
\DeclareFontShape{U}{matha}{m}{n}{
      <5> <6> <7> <8> <9> <10> gen * matha
      <10.95> matha10 <12> <14.4> <17.28> <20.74> <24.88> matha12
      }{}
\DeclareSymbolFont{matha}{U}{matha}{m}{n}

\DeclareMathSymbol{\Lt}{3}{matha}{"CE}
\DeclareMathSymbol{\Gt}{3}{matha}{"CF}

\newcommand{\bbf}{\mathbf{b}}
\newcommand{\pbf}{\mathbf{p}}

\newcommand{\xbf}{\mathbf{x}}

\newcommand{\wbf}{\mathbf{w}}

\newcommand{\Gcal}{\mathcal{G}}

\DeclareMathOperator*{\argmax}{\mathrm{arg\,max}}
\def\BibTeX{{\rm B\kern-.05em{\sc i\kern-.025em b}\kern-.08em
    T\kern-.1667em\lower.7ex\hbox{E}\kern-.125emX}}
\begin{document}

\title{An Optimal Control Approach for Inverse Problems with Deep Learnable Regularizers}
\author[1]{Wanyu Bian}
\affil[1]{Department of Mathematics, University of Florida, Gainesville, FL, USA}
\maketitle

\begin{abstract}
This paper introduces an optimal control framework to address the inverse problem using a learned regularizer, with applications in image reconstruction. We build upon the concept of Learnable Optimization Algorithms (LOA), which combine deep learning  with traditional optimization schemes to improve convergence and stability in image reconstruction tasks such as CT and MRI. Our approach reformulates the inverse problem as a variational model where the regularization term is parameterized by a deep neural network (DNN). By viewing the parameter learning process as an optimal control problem, we leverage Pontryagin's Maximum Principle (PMP) to derive necessary conditions for optimality. We propose the Method of Successive Approximations (MSA) to iteratively solve the control problem, optimizing both the DNN parameters and the reconstructed image. Additionally, we introduce an augmented reverse-state method to enhance memory efficiency without compromising the convergence guarantees of the LOA framework. 
\end{abstract}

\section{Learnable Optimization Algorithms for Solving Inverse Problems}
Deep learning (DL) methods for image processing  have evolved rapidly, offering significant improvements over traditional techniques\cite{zhan2021deepmtl,zhan2022deepmtl,li2024lr,liu2024siambrf,sun2024visual,10594167,ding2024confidence,ni2024earnings,zhou2024reconstruction,yang2024augmentation,ni2024timeseries,ding2024style,ding2024llava,li2024vqa,li2023deception,li2019segmentation}.Deep learning and optimization techniques also applied to complex problems in wireless communication, quantum networks, and sensor networks\cite{zhan2024optimizing,zhan2020efficient,fan2024optimized,ghaderibaneh2024deepalloc}.  One of the most successful DL-based approaches for CT and MRI reconstruction is known as unrolling, which mimics traditional optimization schemes like proximal gradient descent used in variational methods but replaces handcrafted regularization with deep networks. These deep networks are utilized to extract features from images or sinograms, enhancing the reconstruction process. Recently, dual-domain methods have emerged, leveraging complementary information from both image and sinogram domains to further improve reconstruction quality. Despite these advancements, DL-based methods face challenges due to their lack of theoretical interpretation and practical robustness. They are often memory inefficient and prone to overfitting, as they superficially mimic optimization schemes without ensuring convergence and stability guarantees. Additionally, convolutional neural networks (CNNs)  have been successfully applied to sparse-view and low-dose data, projection domain synthesis, post-processing, and prior learning in iterative methods, demonstrating better performance than analytical approaches.

Recently, a new class of DL-based methods known as learnable optimization algorithms (LOA) has been developed for image reconstruction with mathematical justifications and convergence guarantees, showing promising advancements. LOAs are  largely explored in solving inverse problems\cite{chen2021learnable,chen2021variational,zhang2022extra} such as  low-does CT reconstruction\cite{zhang2021provably,zhang2021nonsmooth,} and MRI reconstructions\cite{bian2020deep,bian2021optimization,bian2022optimization,bian2022learnable,bian2022optimal,bian2023magnetic,bian2024diffusion,bian2024multi,bian2024brief,bian2024review,bian2024improving,bian2024review_MDPI}. These methods originate from a variational model where the regularization is parameterized as a deep network with learnable parameters, leading to a potentially nonconvex and nonsmooth objective function. LOAs aim to design an efficient and convergent scheme to minimize this objective function, resulting in a highly structured deep network.  The parameters of this network are inherited from the learnable regularization and are adaptively trained using data while retaining all convergence properties. This approach has been applied to develop sparse-view CT reconstruction and MRI reconstruction methods, using learnable regularizations.  

\section{Solving the Inverse Problem from The Optimal Control Viewpoint}
In this work we propose an optimal control framework to solve the inverse problem by minimize the following varitional model with learned parameters:
\begin{equation}\label{eq:inverse_problem}
  \mathcal{E}(\xbf; \theta) = \mathcal{D}(\xbf, \bbf) + \lambda \mathcal{R} (\xbf; \theta),
\end{equation}
where $\mathcal{D}$ is the data fidelity term to guarantee the recovered signal $\xbf$ is consistent with the measurement $\bbf$, which is usually with the quadratic form $ \frac{1}{2} \| A \xbf - \bbf \|_2^2$ for the linear sampling problem and $A$ is the corresponding measurement matrix.
Here $\mathcal{R}$ is the regularization that aims to incorporate proper prior knowledge of the underlying $\mathbf{x}$, which is realized as a well-structured DNN with parameter $\theta$ to be learned. 

Minimizing the variational model \eqref{eq:inverse_problem} with the gradient flow method \cite{gradient_flow_ref}, we obtain the following ODE
\begin{equation}\label{eq:gradient_flow}
\dot{\xbf}(t) = f ({\xbf}(t), \theta)  := -  \nabla _\xbf \mathcal{E}(\xbf (t) ; \theta), \ \ t \in [0, T].
\end{equation}

We can cast the process of learning the undetermined parameters as an optimal control problem $\mbox{(P)}$ with control parameters $\theta$:
\begin{align*}
\min_{\theta \in  \Theta } \quad & J(\theta) :=  \Phi(\xbf(T)) + \int_{0}^{T} r(\xbf(t),\theta) \dif t, \qquad \mbox{(P)}\\
\mbox{s.t.}\quad & 
\begin{cases}
\dot{\xbf}(t) = f(\xbf(t),\theta), & \quad 0 \le t \le T, \\
\xbf(0) = \xbf^{0}, & 
\end{cases} \qquad \mbox{(ODE)}
\end{align*}
where $\Phi(\xbf(T))$ is the terminal cost and $r(\xbf(t),\theta)$ is the running cost which usually plays role as the regularizer. Besides, $\Theta$ is the admissible set for the control parameters $\theta$ and $\xbf^{0}$ is the initial.

\section{Pontryagin's Maximum Principle}
First we give the definition of the Hamiltonian function $H: \mathbb{R}^n \times \mathbb{R}^n \times \Theta \rightarrow \mathbb{R}$ by
\begin{equation} \label{eq:hamiltonian}
    H(\xbf;\pbf;\theta) = \pbf ^\top f(\xbf,\theta) - r(\xbf,\theta),
\end{equation}
based on which we state the Pontryagin's Maximum Principle (PMP) below
\begin{theorem} [Pontryagin's Maximum Principle, Informal Statement] \label{theorem:pmp}
If $\theta^{*}$ optimizes the optimal control problem \textbf{(P)}, and $\xbf^{*}(t)$ is its corresponding state trajectory. Then, there exists an absolutely continuous co-state process $\pbf^{*}(t), t \in [0, T]$ such that
\begin{equation*}
\label{eq:ode-p}
\begin{cases}
\dot{\xbf}^{*}(t) = \nabla_{\pbf}H(\xbf^{*}(t);\pbf^{*}(t);\theta^{*}), \ \ \xbf^{*}(0) = \xbf_0,
 \qquad \mbox{(ODE)}
\\
\dot{\pbf}^{*}(t) = - \nabla_{\xbf}H(\xbf^{*}(t);\pbf^{*}(t);\theta^{*}),  \pbf^{*}(T) = - \nabla_{\xbf} \Phi(\xbf^{*}(T))^{\top},
\qquad \mbox{(ADJ)}
\\
H(\xbf^{*}(t);\pbf^{*}(t);\theta^{*}) = \max_{\theta} H(\xbf^{*}(t);\pbf^{*}(t);\theta),  \qquad \mbox{(MAX)}
\end{cases}
\end{equation*}
are satisfied.
\end{theorem}

\section{Algorithms of Successive Approximations}
\subsection{Basic MSA}
In light by the PMP, the basic Method of Successive Approximations (MSA) proposed in \cite{JMLR:v18:17-653} are summarized in Algorithm \ref{alg:msa} with an initial guess $\theta^0$ of the control parameter.
\begin{algorithm}[tbh]
\caption{The Basic MSA \cite{JMLR:v18:17-653} to Solve Inverse Problem}
\label{alg:msa}
\begin{algorithmic}[1]
\STATE \textbf{Initialize} $\theta^0$.
\FOR{$k=1,2,\dots,K$ (training iterations)}
\STATE Solve \ $\dot{\xbf}^k(t) = f ({\xbf}^k(t), \theta^{k-1})$, \ \ $\xbf^k(0) = \xbf_{0}$;
\STATE Solve \ $\dot{\pbf}^k(t) = -\nabla_\xbf H ({\xbf}^k(t), \pbf^k(t), \theta^{k-1})$, \ \ $\pbf^k(T) = -\nabla \Phi (\xbf^k(T))^{\top}$;
\STATE Set \ $\theta^{k} = \argmax_{\theta} \int_0^{T} H ({\xbf}^k(t), \pbf^k(t), \theta) dt$;
\ENDFOR
\STATE \textbf{output:} $\theta^{K}$.
\end{algorithmic}
\end{algorithm}

If we apply a gradient ascend step to solve the maximization of the Hamiltonian in Line 5 of Algorithm \ref{alg:msa}, we will have the following update for $\theta$:
\begin{equation} \label{eq:gradient_ascend}
\theta^{k} =  \theta^{k-1} + \eta_k \nabla_{\theta} \int_0^{T}  H ({\xbf}^k(t), \pbf^k(t), \theta^{k-1})dt,
\end{equation}
where $\eta$ is the step size. As proved in \cite{JMLR:v18:17-653}, with a gradient ascent \eqref{eq:gradient_ascend} to maximize the Hamiltonian, the Algorithm \ref{alg:msa} is equivalent to gradient descent with back-propagation (BP).
The same as BP, one bottleneck of the basic MSA is the linear memory cost $\mathcal{O}(T)$ to cache all the intermediate states $\{(\xbf^k(t), \pbf^k(t)): t \in [0, T] \}$. 
\subsection{MSA with Augmented Reverse-State}
As an alternative perspective, \eqref{eq:gradient_ascend} can also be viewed as solving the following ODE backward in time
\begin{equation} \label{eq:ode_backward}
\theta^k = \overline\theta^{k}(0),
\ \  \mbox{where} \ \
\frac{d \overline\theta^{k}(t)}{d t}   =  \eta_k \nabla_{\theta} H ({\xbf}^k(t), \pbf^k(t), \overline\theta^{k}(T)), \ \  \mbox{with} \ \  \overline\theta^{k}(T)  =  \theta^{k-1}.
\end{equation}

Explain:
if we integrate \eqref{eq:ode_backward} backward w.r.t. time $t$, we get
\begin{subequations}
\begin{align}
\overline\theta^{k}(t)  &=  \overline\theta^{k}(T) + \eta_k \int_T^t  \nabla_{\theta} H ({\xbf}^k(t), \pbf^k(t), \overline\theta^{k}(T))dt\\
  &=  \theta^{k-1} + \eta_k \int_T^t  \nabla_{\theta} H ({\xbf}^k(t), \pbf^k(t), \theta^{k-1})dt \qquad \mbox{(Because we set $\overline\theta^{k}(T)  =  \theta^{k-1}$)} \label{eq:prox}
\end{align}
\end{subequations}
Then
\begin{equation}
\theta^k = \overline\theta^{k}(0) 
  =  \theta^{k-1} + \eta_k \int_T^0 \nabla_{\theta} H ({\xbf}^k(t), \pbf^k(t), \theta^{k-1})dt
  \end{equation}
  which is exactly the same as the gradient ascent of the Hamiltonian defined in \eqref{eq:gradient_ascend}.

From the definition of $H$, we know that $\nabla_{\theta} H ({\xbf}, \pbf, \theta) = \pbf^{\top} \nabla_{\theta}f(\xbf, \theta)$.

Instead of solving for co-state $\pbf$ and control $\theta$ separately, here we can solve for the augmented reverse-state $[\pbf, \overline\theta]$ backward in time 
\begin{equation} \label{eq:augmented_reverse_state}
\theta^k = \overline\theta^{k}(0),
\ \  \mbox{where} \ \
\frac{d}{d t} \left[\begin{array}{c} \pbf^k(t)  \\ \overline\theta^{k}(t) \end{array}\right] = \left[\begin{array}{c} -\nabla_\xbf H ({\xbf}^k(t), \pbf^k(t), \overline\theta^{k}(T)) \\ \eta_k \nabla_{\theta} H ({\xbf}^k(t), \pbf^k(t), \overline\theta^{k}(T)) \end{array}\right], \ \ \left[\begin{array}{c} \pbf^k(T)  \\ \overline\theta^{k}(T) \end{array}\right] = \left[\begin{array}{c} -\nabla \Phi (\xbf^k(T))^{\top}  \\ \theta^{k-1} \end{array}\right],
\end{equation} 
which can be easily verified that its numerical result is identical to Algorithm \ref{alg:msa} with Line 5 solved by \eqref{eq:gradient_ascend}. But the merit is that solving \eqref{eq:augmented_reverse_state} only demands constant memory $\mathcal{O}(1)$ to cache the augmented reverse-state $[\pbf, \overline\theta]$. However, in \eqref{eq:augmented_reverse_state} we still need linear memory $\mathcal{O}(T)$ to store all intermediate forward state $\{\xbf^k(t): t \in [0, T] \}$. One way to tackle this issue is to augment it with the reverse state and solve for $\xbf^k(t)$ backward as well with backward initial $\xbf^k(T)$ computed from the forward pass \cite{NEURIPS2018_69386f6b}, as shown in \eqref{eq:augmented_reverse_state_three}.
\begin{equation} \label{eq:augmented_reverse_state_three}
\theta^k = \overline\theta^{k}(0),
\ \  \mbox{where} \ \
\frac{d}{d t} \left[\begin{array}{c} \xbf^k(t) \\ \pbf^k(t)  \\ \overline\theta^{k}(t) \end{array}\right] = \left[\begin{array}{c} \nabla_\pbf H ({\xbf}^k(t), \pbf^k(t), \overline\theta^{k}(T)) \\ -\nabla_\xbf H ({\xbf}^k(t), \pbf^k(t), \overline\theta^{k}(T)) \\ \eta_k \nabla_{\theta} H ({\xbf}^k(t), \pbf^k(t), \overline\theta^{k}(T)) \end{array}\right], \ \ \left[\begin{array}{c} \xbf^k(T) \\ \pbf^k(T)  \\ \overline\theta^{k}(T) \end{array}\right] = \left[\begin{array}{c} \xbf^k(T) \\ -\nabla \Phi (\xbf^k(T))^{\top}  \\ \theta^{k-1} \end{array}\right].
\end{equation} 
The whole process of \eqref{eq:augmented_reverse_state_three} has constant memory cost $\mathcal{O}(1)$ which frees the space of $\{\xbf^k(t): t \in [0, T] \}$. But as a sacrifice we do trade the time for space as solving for \eqref{eq:augmented_reverse_state_three} requires the re-computation of $\xbf^k(t)$. 

If we replace the Line 4-5 in Algorithm \ref{alg:msa} by \eqref{eq:augmented_reverse_state} or \eqref{eq:augmented_reverse_state_three}, we can get an more memory-efficient algorithm with identical numerical result. 

\subsection{MSA with Backward Control Flow}
As a careful observation in \eqref{eq:ode_backward}-\eqref{eq:augmented_reverse_state_three} on the right-hand side of the ODEs, all partial derivatives of $H$ are evalauted on the initial state $\overline\theta^{k}(T)$. In this work, we further free this constraint and give more freedom to the dynamical system. So instead of computing the partial derivative of $H$ on the initial $\overline\theta^{k}(T)$, we compute it on the intermediate state $\overline\theta^{k}(t)$. As discussed in \cite{pmlr-v80-li18b}, the basic idea of the successive approximation methods is to find the optimal parameters from the guess by successive projections onto the manifold defined by the ODEs. So intuitively a better guess will contribute to better convergence performance and even better result. Our modification here gives a better guess to the optimal parameter $\theta^*$  for the Hamiltonian $H$ at the intermediate time $t$, since $\overline\theta^{k}(t)$ is optimized backward in time so $\overline\theta^{k}(t)$ is usually a better estimation point than $\overline\theta^{k}(T)$. Along with \eqref{eq:augmented_reverse_state}, we summarize our proposed algorithm in Algorithm \ref{alg:msa_backward_state}. The forward pass is to compute the trajectory of ${\xbf}^k(t)$, the backward pass is to compute the gradient flow for the co-state ${\pbf}^k(t)$ and control $\overline\theta^{k}(t)$. We would like to point out that here we did not extend \eqref{eq:augmented_reverse_state_three} to the algorithm, because with changing to $\overline\theta^{k}(t)$, the trajectory of $\xbf^k(t)$ might flow to somewhere else, we cannot guarantee $\xbf^k(0) = \xbf_{0}$ anymore. And this might cause the algorithm unstable. We will leave this problem to the future work.
\begin{algorithm}[tbh]
\caption{The MSA with Reverse Augmented State and Control Flow to Solve Inverse Problem}
\label{alg:msa_backward_state}
\begin{algorithmic}[1]
\STATE \textbf{Initialize} $\theta^0$.
\FOR{$k=1,2,\dots,K$ (training iterations)}
\STATE Solve \ $\dot{\xbf}^k(t) = f ({\xbf}^k(t), \theta^{k-1})$, \ \ $\xbf^k(0) = \xbf_{0}$; 
\STATE Solve $\frac{d}{d t} \left[\begin{array}{c} \pbf^k(t)  \\ \overline\theta^{k}(t) \end{array}\right] = \left[\begin{array}{c} -\nabla_\xbf H ({\xbf}^k(t), \pbf^k(t), \overline\theta^{k}(t)) \\ \eta_k \nabla_{\theta} H ({\xbf}^k(t), \pbf^k(t), \overline\theta^{k}(t)) \end{array}\right], \ \ \left[\begin{array}{c} \pbf^k(T)  \\ \overline\theta^{k}(T) \end{array}\right] = \left[\begin{array}{c} -\nabla \Phi (\xbf^k(T))^{\top}  \\ \theta^{k-1} \end{array}\right]$;
\STATE Set $\theta^k = \overline\theta^{k}(0)$.
\ENDFOR
\STATE \textbf{output:} $\theta^{K}$.
\end{algorithmic}
\end{algorithm}

\subsection{Time Discretization}
We discretize $[0, T]$ to be $0, 1, 2 , ..., T$ and employ the explicit Euler method for the forward ODE and Verlet method \cite{Ascher1998ComputerMF} for the backward augmented ODE, then we have the discretized version of Algorithm \ref{alg:msa_backward_state}, which is summarized in Algorithm \ref{alg:discretized}.
\begin{algorithm}[h]
\caption{The Discretized Version of Algorithm \ref{alg:msa_backward_state} to Solve Inverse Problem}
\label{alg:discretized}
\begin{algorithmic}[1]
\STATE \textbf{Initialize} $\theta^0$.
\FOR{$k=1,2,\dots,K$ (training iterations)}
\STATE \textbf{Set} $\xbf^k_0 = \xbf_{0}$;
\FOR{$t=0,1,2,\dots,T-1$}
\STATE $\xbf^{k}_{t+1} = \xbf^k_t + \tau_k f ({\xbf}^k_t, \theta^{k-1})$;
\ENDFOR
\STATE \textbf{Set} $\pbf^k_{T} = -\nabla \Phi (\xbf^k_T)^{\top}$ \textbf{and} $ \overline\theta^{k}_T =  \theta^{k-1}$;
\FOR{$t=T-1,T-2,\dots,0$}
\STATE $\overline\theta^{k}_t = {\overline\theta}^{k}_{t+1} + \eta_k \nabla_\theta H ({\xbf}^k_{t}, \pbf^k_{t+1}, {\overline\theta}^{k}_{t+1})$;
\STATE $\pbf^{k}_t = \pbf^k_{t+1} -\nabla _\xbf H ({\xbf}^k_{t}, \pbf^k_{t+1}, \overline\theta^{k}_{t})$;
\ENDFOR
\STATE \textbf{Set} $\theta^{k+1} = \overline\theta^{k+1}_0$;
\ENDFOR
\STATE \textbf{output:} $\theta^K$.
\end{algorithmic}
\end{algorithm}

\section{Design of the Regularizer in the Variational Model}
In this study, we design the regularizer $\mathcal{R}$ to take the form $\mathcal{R}(\xbf) = \psi(\Gcal(\xbf))\in \Re^{+}$, where $\Gcal: \Re^n \rightarrow \Re^n$ is a feature extraction operator achieved by deep neural network. To enforce the sparsity, the function $\psi$ is defined as $ \psi(\xbf)=\sum_{i=1}^{n} \log \cosh (\xbf_i)$, where $\log \cosh (x)$ is a twice differentiable function approximately equal to $|x| - \log(2)$ for large $x$ and to $x^2/2$ for small $x$.
Throughout this work, we parameterize the feature extraction operator $\Gcal$ as a vanilla $l$-layer convolutional neural network without bias separated by componentwise activation function as follows: 
\begin{equation}\label{eq:g}
  \Gcal(\xbf) = \wbf_l * \sigma ... \ \sigma ( \wbf_3 * \sigma ( \wbf _2 * \sigma ( \wbf _1 * \xbf ))),
\end{equation}
where $ \{\wbf _k \}_{k = 1}^{l}$ denote the convolution weights and $*$ denotes the convolution operation. 
Specifically, we parameterize the first convolution $\wbf _0$ to be $d$ kernels of size $3\times 3$ and the last one $\wbf _l$ with $1$ kernel of size $3\times 3\times d$.
Besides, all hidden layers $ \{\wbf _k \}_{k = 2}^{l-1}$ correspond to convolutions with $d$ kernels of size $3\times 3\times d$.
Here, $\sigma$ represents a componentwise activation function which is twice differentiable. In this work we adopt the twice differentiable function sigmoid-weighted linear unit (SiLU) \cite{elfwing2018sigmoid} as the activation.

\bibliographystyle{IEEEtran}
\bibliography{ref}

\end{document}